\documentclass[twocolumn]{autart}    % Enable this line and disable the 
\newcommand{\0}{{\bf 0}}

\newcommand{\ba}{\begin{array}}
\newcommand{\ea}{\end{array}}
\newcommand{\bc}{\begin{center}}
\newcommand{\ec}{\end{center}}

\newcommand{\beqn}[1]{\begin{equation}\label{#1}}
\newcommand{\eeqn}{\end{equation}}
\newcommand{\be}{\begin{equation}}
\newcommand{\ee}{\end{equation}}

\newcommand{\beqnn}{\begin{eqnarray}}
\newcommand{\eeqnn}{\end{eqnarray}}

\newcommand{\E}{{\rm E}}

\newcommand{\T}{{\rm T}}
\usepackage{graphicx}          % Include this line if your
\usepackage{amsfonts}
\usepackage{amsmath} % assumes amsmath package installed
\usepackage{amssymb}  % assumes amsmath package installed
\usepackage{algorithm}
\usepackage{savesym}
\savesymbol{AND}
\usepackage{algorithmic}
\usepackage{enumerate}
\usepackage{color}
\usepackage{cite,etoolbox}
\usepackage{epstopdf}
\usepackage{booktabs}

\newcommand*{\QEDA}{\hfill\ensuremath{\blacksquare}}

\newcommand{\colt}{\mathtt{col}}
\begin{document}
\begin{frontmatter}
\runtitle{Adaptive model predictive control for a class of constrained linear systems with parametric uncertainties}  % Running title for regular 

\title{Adaptive Model Predictive Control for A Class of Constrained Linear Systems with Parametric Uncertainties\thanksref{footnoteinfo}} % Title, preferably not more 

\thanks[footnoteinfo]{This work was supported by the Natural Sciences and Engineering Research Council of Canada (NSERC). This paper was not presented at any IFAC meeting. Corresponding author Y. Shi. Tel. +1-250-853-3178. Fax +1-250-721-6051.}

\author[Paestum]{Kunwu Zhang}\ead{kunwu@uvic.ca},    % Add the 
\author[Paestum]{Yang Shi}\ead{yshi@uvic.ca}              % e-mail address 
%\author[Baiae]{Publius Maro Vergilius}\ead{vergilius@culture.ir}  % (ead) as shown

\address[Paestum]{Department of Mechanical Engineering, University of Victoria, Victoria, B.C., Canada, V8W 3P6}% Please supply                                              
%\address[Rome]{Senate House, Rome}             % full addresses
%\address[Baiae]{The White House, Baiae}        % here.

                             % keyword list or with the 
 % Five to ten keywords,                                          % help of the Automatica 
                                          % keyword wizard

\begin{abstract}                          % Abstract of not more than 200 words.
This paper investigates adaptive model predictive control (MPC) for a class of constrained linear systems with unknown model parameters. We firstly propose an online strategy for the estimation of unknown parameters and uncertainty sets based on the recursive least square technique. {Then the estimated unknown parameters and uncertainty sets are employed in the construction of homothetic prediction tubes for robust constraint satisfaction. By deriving non-increasing properties on the proposed estimation routine, the resulting tube-based adaptive MPC scheme is recursively feasible under recursive model updates, while providing the less conservative performance compared with the robust tube MPC method.} Furthermore, we theoretically show the perturbed closed-loop system is asymptotically stable under standard assumptions. Finally, numerical simulations and comparisons are given to illustrate the efficacy of the proposed method.
\end{abstract}

\end{frontmatter}
\section{Introduction}\label{sec1:introduction}
Model predictive control (MPC) has become one of the most successful methods for multivariable control systems since it provides an effective and efficient methodology to handle complex and constrained systems \cite{acipl_li_robust_2016}. The main insight of MPC is to obtain a sequence of optimal control actions over the prediction horizon by solving an optimization problem. The prediction employed in MPC is conducted based on an explicit system model. Therefore having an accurate model is critical for {achieving} the desirable performance. However, various categories of uncertainties, such as the measurement noise and the model mismatch, are inevitable in practical control problems. Although standard MPC, which is designed for the nominal system {model} without considering uncertainties, {has certain inherent robustness against} sufficiently small disturbances under certain conditions \cite{yu_inherent_2014Automatica}, its performance may be unacceptable for many practical applications due to the limited robustness, therefore robust MPC has attracted considerable attention in recent years \cite{acipl_Li_Robust_2014,acipl_Liu_Robust_2018,mayne_mpcsurvey_2014auto}. {Yet,} the robust MPC method is generally developed based on the {given bound of uncertainties}, its performance is relatively conservative if the uncertainties are constant or slowly changing. To improve the performance, a general solution is to reformulate the robust MPC scheme manually based on different description of uncertainties \cite{qin_survey_2003CEP}, however, which is relatively resource-intensive and time-consuming. Alternatively, a promising solution is to allow for {online model adaptation in the MPC framework}, which is termed as adaptive MPC.  

In recent years, adaptive MPC has drawn increasing attention since it provides a promising solution to reduce conservatism of robust MPC by {incorporating} system identification into the robust MPC framework. Mayne and Michalska firstly proposed an adaptive MPC method in \cite{mayne_adaptive1993CDC} for input-constrained nonlinear uncertain systems, where the convergence of parameter estimates can be guaranteed if the MPC problem is recursively feasible. Later in \cite{fukushima_Adaptive_2007automatica}, a data selection mechanism is considered to improve the convergence performance of parameter estimates for linear systems, then the estimated error bound is employed to construct the comparison model used in the robust MPC framework. But this method relies on the system {model} represented in a controllable canonical form. An alternative of {fulfilling} the persistent excitation (PE) condition is to impose an additional constraint on system states \cite{marafioti_persistently_2014IJACSP} or control inputs \cite{vicente_stabilizing_2019SCL} to the MPC optimization problem. With this strategy, the convergence of parameter estimates can be ensured, but the system state can only be stabilized in a small region around the origin due to the presence of constraints from the PE condition. 

In the literature, another category of research on adaptive MPC is to combine set-membership identification with robust MPC. In \cite{adetola_adaptiveMPC_2009SCL}, an ellipsoidal uncertainty set is constructed based on the recursive least square (RLS) technique, then a stabilizing min-max MPC scheme is developed for constrained continuous-time nonlinear systems. The discrete-time version of this approach is presented in \cite{adetola_robust_2011IJACSP}. The polytope based set-membership identification is considered in \cite{tanaskovic_adaptive_2014Automatica}, where an adaptive output feedback MPC {approach} is designed for constrained stable finite impulse response (FIR) systems. This method has been extended to handle chance constraints \cite{Bujarbaruah_Adaptive_2018ACC} and time-varying systems \cite{tanaskovic_adaptive_2018ECC}. A combination of the set-membership identification and homothetic tube MPC is proposed in \cite{lorenzen_adaptive_2017IFAC,lorenzen_robust_2019Automatica}, in which the worst-case realization of the uncertainty is considered based on a set-based state prediction and uncertainty estimation. Recently, incorporating machine learning techniques with robust MPC has also attracted much attention. The true model is described by a nominal model plus a learned model, e.g., the Gaussian process \cite{ostafew_learning_2016JFR} and the neural network \cite{nagabandi_neural_2017arXiv,zhang_safe_2019arxiv}. Although aforementioned works on machine learning based methods have showed empirical success, how {to theoretically} guarantee the closed-loop stability and recursive feasibility with desired estimation performance is still a major challenge. {The application of adaptive MPC to repetitive or iterative processes can be found in \cite{bujarbaruah_adaptive_2018,rosolia_learning_2017TAC}.}

In this work, we propose a computationally tractable adaptive MPC algorithm for a class of constrained linear systems subject to parametric uncertainties. Similar to \cite{adetola_adaptiveMPC_2009SCL}, the proposed method uses an RLS based estimator to identify the unknown system parameters. {Note that the estimated uncertainty set in \cite{adetola_adaptiveMPC_2009SCL} is employed to update the min-max optimization problem for robust constraint satisfaction, which is non-convex and computationally complicated.} Alternatively, the proposed work employs the tube MPC technique ,e.g., \cite{fleming_robust_2015TAC,langson_robust_2004Automatica,hanema_stabilizing_2017automatica,rakovic_homothetic_2012automatica}, to handle the uncertainty, which has a comparable computational complexity to standard MPC. Recently, there are some novel adaptive MPC strategies \cite{lorenzen_adaptive_2017IFAC,lorenzen_robust_2019Automatica} combining the homothetic tube MPC technique, e.g., \cite{rakovic_homothetic_2012automatica}, with the set-membership identification, where the sequence of state tubes $\{\mathcal{X}_{l|k}\}$ is developed with the form $\mathcal{X}_{l|k} = z_{l|k}+\sigma_{l|k}\mathcal{X}_0$ to guarantee the robust constraint satisfaction.  Here, $z_{l|k}$ is the nominal system state, $\mathcal{X}_0$ is a given set and $\sigma_{l|k}$ is a scalar to be optimized by the MPC optimization problem. {It can be seen that the tube cross sections are shaped by the set $\mathcal{X}_0$, translated and scaled by the MPC optimization problem.} The set $\mathcal{X}_0$ is calculated offline according to the initial knowledge of the uncertainty set, {which may be conservative under recursive updates of the uncertainty set. Inspired by the tube MPC approach in \cite{fleming_robust_2015TAC}, we construct the homothetic tubes in this work, where both the size and shape of the tube cross sections are optimized via the MPC optimization problem.} Consequently, it will {promisingly} lead to control performance improvement by using the proposed method. The main contribution of this work is to extend the robust MPC framework in \cite{fleming_robust_2015TAC} to allow for {online model adaptation}, while guaranteeing the closed-loop stability and recursive feasibility. Compared with the methods in \cite{lorenzen_adaptive_2017IFAC,lorenzen_robust_2019Automatica}, the proposed {approach} introduces additional decision variables in the MPC optimization problem to optimize both the shape and size of the tube cross sections, {resulting in the reduced conservatism.} 
In addition, to provide a trade-off between the computational complexity and conservatism, a specialization of {the} proposed adaptive method is also given with reduced computational complexity and {comparable control performance}. A numerical example and comparison study are given to illustrate the benefits of the proposed method.
 
The remainder of this paper is organized as follows: Section \ref{sec1:formulation} 
demonstrates the problem formulation. In Section \ref{sec1:estimation}, the
estimation of the unknown parameter and the uncertainty set are discussed. An adaptive MPC algorithm is presented in Section \ref{sec1:LCMPC}, followed by the analysis of closed-loop stability and recursive feasibility. Simulation and comparison studies are illustrated in Section \ref{sec1:simulation}. Finally Section \ref{sec1:conclusion} concludes this work. 
\section{Problem Formulation}\label{sec1:formulation}
\subsection{Notation}
Let $\mathbb{R}, \mathbb{R}^{n}$ and $\mathbb{R}^{m\times n}$ denote the sets of real numbers, column real vectors with $n$ components and real matrices consisting of $n$ columns and $m$ rows, respectively.  The notation $\mathbb{N}$ denotes the set of non-negative integers, and  $\mathbb{N}_a^b = \{x\in\mathbb{N}|a\leq x\leq b\}$. Given a vector $x\in\mathbb{R}^{n}$, the Euclidean norm and infinity norm of $x$  are denoted by $\|x\|$ and $\|x\|_\infty$, respectively. We define $\|x\|_Q = x^\T Q x$. The Pontryagin difference of sets $X\subseteq \mathbb{R}^n$ and $Y\subseteq \mathbb{R}^n$ is denoted by $X \ominus Y = \{z\in\mathbb{R}^n: z + y\in X; \forall y \in Y \}$, and the Minkowski sum is $X\oplus Y = \{x+y|x\in X, y\in Y \}$. {The column operation $\colt(\cdots)$ is defined as $\colt(x_1,x_2,\cdots,x_n) = [x_1^\T,x_2^\T,\cdots,x_n^\T]^\T$.} We use $I_{n}$ to denote an identity matrix of size $n$. For an unknown vector $\theta$, the notations $\hat{\theta}$ and $\theta^*$ represent its estimation and real value, respectively. Then the estimation error is defined as $\tilde{\theta} = \theta^*-\hat{\theta}$.
\subsection{Problem setup}
Consider a discrete-time linear time-invariant (LTI) system with an unknown parameter $\theta\in\mathbb{R}^{n_\theta}$	 
\begin{equation}\label{sec1:sys}
	x_{k+1}=A(\theta)x_k+B(\theta)u_k,
\end{equation}
subject to a mixed constraint
\begin{equation}\label{sec1:mix_constaints} 
	\mathcal{M} = \{({x_k,u_k})|Fx_k+Gu_k\leq \mathbf{1}\},  	
\end{equation}
where $x_k\in\mathbb{R}^{n_x}$ and $u_k\in\mathbb{R}^{n_u}$ are the system state and input, respectively. The matrices $A(\theta)$ and $B(\theta)$ are the real affine functions of  $\theta$, i.e., $A(\theta)=A_0+ \sum_{i=1}^{n_\theta}A_i\theta_i, B(\theta)=B_0+\sum_{i=1}^{n_\theta}B_i\theta_i$. $\theta = \colt(\theta_1,\theta_2,\cdots,\theta_{n_\theta})$ is the vector of unknown parameters, which is assumed to be uniquely identifiable \cite{saccomani_parameter_2003}. {It is assumed that the parameter $\theta$ is bounded by a given set $\Theta_0 = \{\theta| \|\theta\|\leq r_0\}$ which contains the real parameter $\theta^*$.} 

In this paper, the goal is to design a state feedback control law for the perturbed and constrained system in (\ref{sec1:sys}) while ensuring the desirable closed-loop {performance} and robust constraint satisfaction by means of adaptive MPC. In particular,
we consider the following parameterization of the control input 
{\begin{equation}\label{sec1:u_dual}
u_k = Kx_k+v_k,
\end{equation}}where $v_k\in\mathbb{R}^{n_u}$ is the decision variable of the MPC optimization problem; $K\in\mathbb{R}^{n_u\times n_x}$ is a prestabilizing state feedback gain such that $\phi(\theta)=A(\theta)+B(\theta)K$ is quadratically stable for all $\theta \in \Theta_0$. 
{
\begin{defn}[\cite{blanchini_set_2008Springer}]
	A polyhedral set $\mathcal{Z}$ is robustly positively invariant (RPI) for the system in (\ref{sec1:sys}) with respect to the constraint (\ref{sec1:mix_constaints}) and the feedback control law $u_k = Kx_k$ iff $(x_k,Kx_k)\in\mathcal{M}$ and $\phi(\theta)x_k\in\mathcal{Z}$ for all $x_k\in\mathcal{Z}$ and $\theta\in\Theta_0$.                                       
\end{defn}
Suppose that $\mathcal{Z}$ is an RPI set for the system in (\ref{sec1:sys}) with respect to the constraint (\ref{sec1:mix_constaints}) and the control law $u_k = Kx_k$, if  $\mathcal{Z}$ contains every RPI set, then $\mathcal{Z}$ is the maximal RPI (MRPI) set for {the} system in (\ref{sec1:sys}). As shown in \cite{pluymers_efficient_2005ACC}, if the MRPI set $\mathcal{Z}$ for {the} system in (\ref{sec1:sys}) exists, it is unique. An example of calculating the MPRI set $\mathcal{Z}$ can be found in \cite{pluymers_efficient_2005ACC}.}

%%%%%%%%%%%%%%%%%%%%%%%%%%%%%%%%%%%%%%%%%%%%%%%%%%%%%%%%%%%%%%%%%%%%%%%%%%%%%%%%
\section{Uncertainty Estimation}\label{sec1:estimation}
In this section, we introduce an online parameter estimation scheme based on the RLS technique with guaranteed non-increasing estimation errors. {Thereafter, in order to reduce conservatism in robust MPC, the approximation of feasible solution set (FSS) of the unknown parameters is presented. Finally, we conclude this section by analyzing the performance of the proposed estimation scheme.}

\subsection{Parameter estimation}\label{sec1:para_est}
Let $g(x_k,u_k)\theta =\sum_{i=1}^{n_\theta}(A_ix_k+B_iu_k)\theta_i$, then we can formulate a regressor model $y_k = g(x_k,u_k)\theta^*$ with $y_k = x_{k+1}-A_0x_k-B_0u_k$ to estimate $\theta^*$ by using the standard RLS method. But the convergence of this solution relies on the PE {condition} of $g(x_k,u_k)$, which cannot be guaranteed if $x_k = 0$ and $u_k = 0$. Similar to \cite{adetola_robust_2011IJACSP}, we introduce the following filter $w_k$ for the regressor $g(x_k,u_k)$ to improve the convergence performance,   
\begin{equation}\label{sec1:est_w}
	w_{k+1}=g(x_k,u_k)-K_ew_k, 
\end{equation}
where $w_0 = 0$ and $K_e$ is a Schur stable gain matrix. Let $\hat{x}_k$ denote the system state estimated at time $k-1$, based on (\ref{sec1:sys}) and (\ref{sec1:est_w}), a state estimator at time $k$ is designed as follows:
\begin{equation}\label{sec1:state_est}
	\begin{array}{l}
		\hat{x}_{k+1}=A_0x_k+B_0u_k+g(x_k,u_k)\hat{\theta}_{k+1}\\
		~~~~~~~~~~~+K_e\tilde{x}_k+K_ew_k(\hat{\theta}_k-\hat{\theta}_{k+1}),
	\end{array}
\end{equation}
where $\tilde{x}_k = x_k-\hat{x}_k$ is the state estimation error. Then subtracting (\ref{sec1:sys}) from (\ref{sec1:state_est}) yields 
\begin{equation}\label{sec1:err_evol}
	\tilde{x}_{k+1}=g(x_k,u_k)\tilde{\theta}_{k+1}-K_e\tilde{x}_k-K_ew_k(\hat{\theta}_k-\hat{\theta}_{k+1}).
\end{equation} In order to establish an implicit regression model for $\hat{\theta}$, we introduce an auxiliary variable $\eta_k$ in the following
\begin{equation}\label{sec1:eta_def}
	\eta_k=\tilde{x}_k-w_k\tilde{\theta}_k.
\end{equation} 
Then by substituting (\ref{sec1:est_w})-(\ref{sec1:err_evol}) into (\ref{sec1:eta_def}), one gets
\begin{equation}\label{sec1:eta_evol}
	\eta_{k+1}=-K_e\eta_k.
\end{equation}
Based on this implicit regression model, we develop the following parameter estimator by using the standard RLS algorithm \cite{johnstone_exponential_1982SCL}
\begin{subequations}\label{sec1:theta_evol}
	\begin{align}
		&\hat{\theta}_{k+1} = \hat{\theta}_k+\Gamma_{k+1}^{-1}w_k^\T(\tilde{x}_k-\eta_k),\label{sec1:theta_evol_a}\\
		&\Gamma_{k+1}=\lambda\Gamma_k+w_k^\T w_k,\label{sec1:theta_evol_Gamma}
	\end{align}
\end{subequations}
where $\Gamma_0=\beta I_{n_\theta}$; $\beta$ is the positive scalar, and $\lambda\in(0,1)$ is the forgetting factor. Then it follows from \cite{johnstone_exponential_1982SCL} that the non-increasing estimation error is guaranteed, and the convergence of parameter estimates $\hat{\theta}_k$ can be achieved if the sequence $w_k$ is persistently exciting. 

By using the proposed estimation mechanism (\ref{sec1:theta_evol}), the convergence of the estimation error $\tilde{\theta}_k$ relies on the persistently exciting sequence of $w_k$ instead of $g(x_k,u_k)$. {Suppose that the system is stable when $k\geq t_s,t_s\in\mathbb{N}_0^\infty,$ and $w_{t_s}\neq0$. According to (\ref{sec1:est_w}), we have $w_{k+1} = - K_ew_k$ for all $k\geq t_s$. Let $\mathbf{w}_k = \{w_{k},w_{k+1},\cdots,w_{k+N_p-1} \}$ with $N_p\in\mathbb{N}_0^\infty$. Then it can be derived that $\mathbf{w}_k\mathbf{w}_k^\T = \mathcal{K}_ew_kw_k^\T\mathcal{K}_e^\T$ for $k\geq t_s$, where $\mathcal{K}_e =\colt( I,-K_e,(-K_e)^2,\cdots,(-K_e)^{N_p-1})$. Since $K_e$ is Schur stable, it is possible to find $N_p, l_p\in\mathbb{N}_0^\infty,t_d\in\mathbb{N}_{t_s}^\infty, \rho_0>0$ and $\rho_1>0$ such that $\rho_1 I_{N_pN_x}>\sum_{j=0}^{l_p-1}(\mathbf{w}_{k+j}\mathbf{w}_{k+j}^\T)>\rho_0I_{N_pN_x}$ for all $k\in\mathbb{N}_{t_s}^{t_d}$. Therefore, the sequence $\mathbf{w}_k$ satisfies the PE condition during a certain period when the system is stable. In addition, it can be derived from (\ref{sec1:theta_evol}) that $\hat{\theta}_{k+1}\approx \hat{\theta}_{k}, \Gamma_{k+1}\approx \lambda\Gamma_{k}$ and the corresponding $\hat{\Theta}_{k}\approx\hat{\Theta}_{k+1}$ when $w_k$ is sufficiently small. {Since $w_k$ is decreasing when the system in (\ref{sec1:sys}) is stable, $\hat{\Theta}_k$ will converge to a fixed set in finite time}.} 
%The proposed estimation mechanism (\ref{sec1:theta_evol}) can provide a better (at least equivalent) convergence {performance} than the standard RLS method. 

\subsection{Uncertainty set estimation}
To bound the unknown parameters, we introduce the following ellipsoidal uncertainty set
\begin{equation}\label{sec1:uncer_set_hat}
\hat{\Theta}_{k} = \{\theta| \|\theta-\hat{\theta}_{k}\|_{\Gamma_{k}}\leq \mathcal{V}_{k} \}.
\end{equation}
where $\mathcal{V}_{k}>0$ is the bound of the estimation error. According to (\ref{sec1:theta_evol_Gamma}), we define the propagation of $\mathcal{V}_{k}$ as $\mathcal{V}_{{k+1}}=\lambda \mathcal{V}_{k}$ with $\mathcal{V}_{{0}}=\bar{\Lambda}(\Gamma_0)r_0^2$, where $\bar{\Lambda}(\Gamma_0)$ is the maximal eigenvalue of $\Gamma_0$.

Let $\Theta_{k}$ denote the FSS of unknown parameters. Since unknown parameters are uniquely identifiable and stay in the \textit{a priori} known set $\Theta_0$, $\Theta_{k}$ must be the subset of $\Theta_0$. Therefore, for all $k\geq 1$, $\Theta_{k}$ is computed as follows
\begin{equation}\label{sec1:uncer_set}
\Theta_{k} = \Theta_{k-1}\cap\hat{\Theta}_{k}.
\end{equation}
By choosing suitable $\hat{\theta}_0,\Gamma_0$ and $\mathcal{V}_{0}$, $\hat{\Theta}_0$ can be equivalent to $\Theta_{0}$. 
The following lemma shows the performance of uncertainty set estimation. 

\begin{lem}\label{sec1:lemma:uncer_set}
	Let $\Theta_{k}$ denote the estimated uncertainty set updated by following (\ref{sec1:est_w})-(\ref{sec1:uncer_set}) at each time instant. Suppose that $\theta^*\in\Theta_0$, then we have $\theta^*\in\Theta_k$ for all $k\geq 0$.
\end{lem}
\begin{pf*}{Proof}
	To prove this lemma, we firstly show that $\theta^*\in\hat{\Theta}_k$ for all $k\geq0$. Let $\mathcal{V}(\tilde{\theta}_k)=\tilde{\theta}_k^\T\Gamma_{k}\tilde{\theta}_k$, then it follows from \cite{johnstone_exponential_1982SCL} that $\mathcal{V}(\tilde{\theta}_k)$ is non-increasing and $\mathcal{V}(\tilde{\theta}_{k})\leq\lambda \mathcal{V}(\tilde{\theta}_{k-1})$. {When $k=0$, the condition $\mathcal{V}(\tilde{\theta}_0) = \tilde{\theta}^\T_0\Gamma_0\tilde{\theta}_0\leq \mathcal{V}_{0}$ holds by using $\|\tilde{\theta}_0\|\leq r_0$. When $k>0$, we still have $\mathcal{V}_{k}\geq \mathcal{V}(\tilde{\theta}_k)$ since $\mathcal{V}_{k} = \lambda^k\mathcal{V}_{{0}}$ and $\mathcal{V}(\tilde{\theta}_{k})\leq\lambda^k\mathcal{V}(\tilde{\theta}_0)$. Therefore, one gets $\mathcal{V}(\tilde{\theta}_{k})\leq\mathcal{V}_{k}$ for all $k\geq 0$. Then according to (\ref{sec1:uncer_set_hat}), it can be derived that $\theta^*\in\hat{\Theta}_{k}$ for all $k\geq 0 $.} Suppose that $\theta^*\in\Theta_{k}$. At next time instant, we have $\theta^*\in\hat{\Theta}_{k+1}$, which implies that $\theta^*\in\Theta_{k}\cap\hat{\Theta}_{k+1} = \Theta_{k+1}$. Hence, it can be concluded that $\theta^*\in\Theta_{k}$ for all $k\geq 0$ if $\theta^*\in\Theta_{0}$.
\QEDA\end{pf*}
	Generally, the tightened state constraints are widely employed in robust MPC to guarantee recursive feasibility and closed-loop stability. These constraints are designed based the {given bounds of uncertainties}. Hence, having an accurate description on the uncertainty is crucial to obtain the desirable closed-loop performance. By incorporating the proposed parameter estimator, it is possible to use the estimated parameters and uncertainty sets at each time instant to obtain more accurate predictions and less conservative tightened state constraints in {robust} MPC, and thus improving the control performance. In the following section, a computationally tractable integration of tube MPC and the proposed estimator is presented.   
%%%%%%%%%%%%%%%%%%%%%%%%%%%%%%%%%%%%%%%%%%%%%%%%%%%%%%%%%%%%%%%%%%%%%%%%%%%%%%%%

\section{Adaptive Model Predictive Control}\label{sec1:LCMPC}
In this section, we present a computationally tractable adaptive MPC algorithm based on the homothetic tube MPC technique. Let $x_{l|k}$ denote the predicted real system state $l$ steps ahead from time $k$ and $x_{l|k} = z_{l|k} + e_{l|k}$, where $z_{l|k}$ and $e_{l|k}$ are the predicted nominal system state and the error state, respectively. Our objective is to design a sequence of state tubes $\{\mathcal{X}_{l|k}\}$ {for robust constraint satisfaction, i.e., the following conditions hold for some $u_{l|k}$:}
\begin{subequations}\label{sec1:tube_cons}
	\begin{align}
	&x_k\in\mathcal{X}_{0|k}\\
	&A(\theta)x+B(\theta)u_{l|k}\in\mathcal{X}_{l+1|k},~\forall x\in\mathcal{X}_{l|k},\theta\in\Theta_{k+1}\\
	&(x,u_{l|k})\in\mathcal{M},\forall x\in\mathcal{X}_{l|k}
	\end{align}
\end{subequations}
Instead of designing the state tube $\mathcal{X}_{l|k}$ directly, in this work we construct the tube cross section $S_{l|k}$ for the error state $e_{l|k}$. {Therefore, the state tube can be established indirectly as} $\mathcal{X}_{l|k} = z_{l|k}\oplus S_{l|k}$. In the following, we present how to design the homothetic tubes according to the estimation of uncertainties.

\subsection{Error tube and constraint satisfaction}
{As mentioned in Section \ref{sec1:para_est}, we predict $\hat{\theta}_{k+1}$ and ${\Theta}_k$ at time $k$ based on the state estimation error $\tilde{x}_k$. Hence the system matrices $A({\hat{\theta}}_{k+1})$ and $B({\hat{\theta}}_{k+1})$ are considered in the following for predicting the nominal system state at time $k$:}
\begin{equation}\label{sec1:sys_nol}
	z_{l+1|k} = {A}_{k+1}z_{l|k} +{B}_{k+1}u_{l|k} 
\end{equation} 
where ${A}_{k+1} = A({\hat{\theta}}_{k+1})$ and ${B}_{k+1} = B({\hat{\theta}}_{k+1})$; $N$ is the prediction horizon and $l\in\mathbb{N}_{0}^{N-1}$.

Then subtracting (\ref{sec1:sys}) from (\ref{sec1:sys_nol}) results in 
\begin{equation}\label{sec1:err_dyn}
	\begin{array}{ll}
		e_{l+1|k}&=x_{l+1|k}-z_{l+1|k}\\
		%&= \phi^*x_{l|k}-\phi_kz_{l|k}+B(\theta^*)v_{l|k}-B_kv_{l|k}\\
		&=\phi^* e_{l|k}+\Delta\phi_{k+1}z_{l|k}+\Delta B_{k+1}v_{l|k},
	\end{array}
\end{equation}
where $\phi^* = A(\theta^*)+B(\theta^*)K, \phi_{k+1} = A_{k+1}+B_{k+1}K$, $\Delta\phi_{k+1} = \phi^*-\phi_{k+1}$ and $\Delta B_{k+1} =B(\theta^*)-B_{k+1}$. {Since $\Theta_k$ is compact and convex, we can find a polytope to over approximate  $\Theta_{k}$ by following the algorithm in \cite{sapatnekar_feasible_1993ISCAS}. Let $\bar{\Theta}_k$ denote the polytopic over approximation of $\Theta_k$, and $\mathtt{Pol}(\cdot)$ is the polytopic approximation operator from the algorithm in \cite{sapatnekar_feasible_1993ISCAS}. Hence $\bar{\Theta}_k$ can be directly calculated as $\bar{\Theta}_k = \mathtt{Pol}(\Theta_{k})$. Due to the recursive set intersection in (\ref{sec1:uncer_set}), we calculate $\bar{\Theta}_k$ indirectly to reduce the computational load, i.e., $\bar{\Theta}_{k} =  \mathtt{Pol}(\hat{\Theta}_{k})\cap\bar{\Theta}_{k-1}$	with $\bar{\Theta}_{0} = \mathtt{Pol}(\Theta_{0})$.} Suppose that $\bar{\Theta}_{k}$ can be equivalently represented by a convex hull  $Co({\hat{\theta}}^j_k)$ where $j\in\mathbb{N}_0^{{n_c}}$ and $n_c$ is an integer denoting the number of extreme points in the convex hull. Hence, a set for the system pair $(A(\theta), B(\theta))$ at time $k$ can be approximated by using a convex hull $Co({A}^j_k, {B}_k^j)$ where ${A}^j_k = A(\hat{\theta}^j_k)$ and  ${B}_k^j = B(\hat{\theta}^j_k)$.  
%Similarly, let  ${\phi}^j_k ={A}^j_k+{B}^j_kK$, then the convex hulls $Co({\phi}_k^j)$ and $Co({B}_k^j)$ that approximate the set of $\phi(\theta)$ and $B(\theta)$ can be found, respectively. 

Inspired by the previous work \cite{fleming_robust_2015TAC}, we consider a polytopic tube with the form $S_{l|k} = \{e_{l|k}|Ve_{l|k}\leq \alpha_{l|k}\}$ for the error $e_{l|k}$ to handle multiplicative uncertainties, where $V\in\mathbb{R}^{n_v\times n_x}$ is a matrix describing the shape of $S_{l|k}$; $\alpha_{k}\in\mathbb{R}^{n_v\times 1}$ is the tube parameter to be optimized. The following proposition shows a sufficient condition for the robust satisfaction of constraint (\ref{sec1:mix_constaints}).   

\begin{prop}\label{sec1:tube_proposition}
	Let $S_{l|k} = \{e_{l|k}|Ve_{l|k}\leq \alpha_{l|k} \}$. Suppose that $e_{l|k}\in S_{l|k}$, then $e_{l+1|k}\in S_{l+1|k}$. In addition, the constraint (\ref{sec1:mix_constaints}) is satisfied at each time instant if the following conditions hold:
	\begin{subequations}\label{sec1:tube}
		\begin{align}
			&\mathbf{1}\geq\begin{cases}\label{sec1:tube1}
				H\alpha_{l|k}+(F+GK)z_{l|k}+Gv_{l|k},l\in\mathbb{N}_0^{N-1}\\
				H\alpha_{l|k}+(F+GK)z_{l|k},l\in\mathbb{N}_N^{\infty}
			\end{cases}\\
			&\alpha_{l+1|k}\geq	H_{k+1}^j\alpha_{l|k}+V(\Delta\phi_{k+1}^jz_{l|k}+\Delta B_{k+1}^jv_{l|k})
			\nonumber\\
			&~~~~~~~~~~~~~~~~~~~~~~~~~~~~~~~~~~~~~~~~~~~l\in\mathbb{N}_0^{\infty},j\in\mathbb{N}_0^{n_c}\label{sec1:tube2}
		\end{align}
	\end{subequations}
	where ${\phi}^j_k ={A}^j_k+{B}^j_kK, \Delta\phi_{k+1}^j = {\phi}_{k+1}^j -\phi_{k+1}$ and $\Delta B_{k+1}^j = {B}_{k+1}^j -B_{k+1}$; $H$ and $H_{k+1}^j$ are non-negative matrices satisfying the conditions
	$HV=F+GK$ and $H_{k+1}^jV=V{\phi}_{k+1}^j$.
\end{prop}  

\begin{pf*}{Proof}
	Consider the uncertain input matrix $B({\theta})$ in the system (\ref{sec1:sys}), this proof is completed by following the proof of Proposition 2 in \cite{fleming_robust_2015TAC}.
\QEDA\end{pf*}

Proposition \ref{sec1:tube_proposition} shows a sequence of tightened sets for the nominal system state. By considering tube parameters $\{\alpha_{l|k}\}$ as extra decision variables of the MPC optimization problem, we can obtain the optimal tube cross sections online. 

According to the proposed parameter estimator, we can obtain the new estimation of the real system with non-increasing estimation error at each time instant. Hence, a time-varying nominal system is used to improve the accuracy of prediction. However, the system is considered to be invariant during the prediction. In order to improve the control performance, a time-varying terminal set is constructed  based on the new estimation of uncertainty, which will be presented in the following.

\subsection{Construction of terminal sets}
Based on Proposition \ref{sec1:tube_proposition}, we define the following dynamics of $z_{l|k}$ and $\alpha_{l|k}$ for $l\in\mathbb{N}_N^\infty$ at time $k$
\begin{subequations}\label{sec1:sys_mode2}
	\begin{align}
		&\alpha_{l+1|k} = \max_{j\in\mathbb{N}_0^{n_c}} \{H_k^j\alpha_{l|k}+V\Delta\phi_{k+1}^jz_{l|k}\}\label{sec1:alpha_dyn},\\
		&z_{l+1|k} = \phi_{k+1}z_{l|k},\label{sec1:z_dyn}
	\end{align}
\end{subequations}
where the maximization is taken for each element in the vector. Let $\mathcal{Z}_k$ denote the polytopic RPI set for the system $x_{k+1} = (A(\theta)+B(\theta)K) x_{k}$ with respect to the uncertainty set $\Theta_{k+1}$. Since $\hat{\theta}_{k+1}\in\Theta_{k+1}$, $\mathcal{Z}_{k}$ is also RPI for the system in (\ref{sec1:z_dyn}).

Define $\mathcal{Z}_{l+1|k}^j$ as $\mathcal{Z}_{l+1|k}^j = \phi_{k+1}^j\mathcal{Z}_{l|k}^j$ with $\mathcal{Z}_{0|k}^j = \mathcal{Z}_k$ for all  $j\in\mathbb{N}_0^{n_c}$, then we have $\mathcal{Z}_{l+1|k}^j\subseteq \mathcal{Z}_{l|k}^j\subseteq\mathcal{Z}_{k}$ since $\phi^j_{k+1}$ is Schur stable for all $j\in\mathbb{N}_0^{n_c}$. Inspired by Proposition 3 in \cite{fleming_robust_2015TAC}, the following proposition is given to construct the invariant set for the system in (\ref{sec1:alpha_dyn}).     
\begin{prop}\label{sec1:bar_alpha}
	Define
	\begin{equation}\label{sec1:fc_prop}
	\begin{array}{l}
		\bar{f}_{l|k}^j = \max_{z\in\mathcal{Z}_{l|k}^j}\{(F+GK)z \},\\\bar{c}_{l|k}^j=\max_{z\in\mathcal{Z}_{l|k}^j}\{V({\phi}^j_{k+1}-\phi_{k+1})z\},\\
		\bar{g}_{l|k}^j=\max_{z_1,z_2\in\mathcal{Z}_{l|k}^j}\{V{\phi}^j_{k+1}(z_1-z_2)\}.
	\end{array}
	\end{equation}
	%	and $H_{l|k}^j$ is nonnegative matrix such that $H^jV=V{\phi}^j_0$.
	The set $\mathcal{A}_k=\{\alpha| \|\alpha\|_\infty\leq\gamma_k, \alpha\geq 0\label{key} \}$ is invariant for the system in (\ref{sec1:alpha_dyn}) while the constraint $H\alpha+(F+GK)z\leq\mathbf{1}$ is satisfied if the following condition holds
	\begin{equation}\label{sec1:gam_bd}
		\bar{\gamma}_{l|k}\geq\gamma_k\geq\underline{\gamma}_{l|k}
	\end{equation}
	where 
	\begin{equation*}
	\begin{array}{l}
			\underline{\gamma}_{l|k} =  \frac{\max_{j\in\mathbb{N}_0^{n_c}} \|\bar{c}_{l|k}^j\|_\infty+\|\bar{g}_{l|k}^j\|_\infty }{1-\max_{j\in\mathbb{N}_0^{n_c}}\|H^j_{k+1}\|_\infty},
	\bar{\gamma}_{l|k} = \frac{1-\max_{j\in\mathbb{N}_0^{n_c}}\|\bar{f}_{l|k}\|_\infty}{\|H\|_\infty}.
	\end{array}
	\end{equation*}
	In addition, there exists a $\gamma_k$ satisfying the condition (\ref{sec1:gam_bd}) if $l$ is sufficiently large.
\end{prop}  
\begin{pf*}{Proof}
	This proposition can be proved by following the proof of Proposition 3 in \cite{fleming_robust_2015TAC}.
\QEDA\end{pf*}

As shown in \cite{fleming_robust_2015TAC}, the invariant set $\mathcal{A}_k$ for the system in (\ref{sec1:alpha_dyn}) is nonempty if 
$\|H_k^j\|_\infty<1$ for all $k\geq 0 $. This condition can be satisfied by choosing the appropriate $V$ such that the set $\{x|Vx\leq \mathbf{1}\}$ is a $\lambda$-contractive set for the system $z_{k+1} = \phi(\theta)z_k,\forall \theta\in\Theta_{0}$. An example of computing {the} matrix $V$ can be found  in \cite{blanchini_set_2008Springer}.

\begin{lem}\label{sec1:horizon_mode2}
	{Given uncertainty sets $\Theta_{k}$ and $\Theta_{k+1}$, assume that the sets $\mathcal{A}_k$ and $\mathcal{A}_{k+1}$ are not empty. Let $M_k$ denote the minimum $l$ that the condition (\ref{sec1:gam_bd}) is satisfied.} Then we have $M_k\geq M_{k+1}$ if the condition $\Theta_{k+1}\subseteq \Theta_{k}$ holds.
\end{lem}
\begin{pf*}{Proof}
	{According to (\ref{sec1:tube}) and (\ref{sec1:fc_prop})}, if $\Theta_{k+1}\subseteq \Theta_{k}$, it can be derived that
	\begin{equation}\label{sec1:gamma_lemma1_1}
		\bar{\gamma}_{l|k+1}\geq \bar{\gamma}_{l|k}, \underline{\gamma}_{l|k+1}\leq \underline{\gamma}_{l|k}.
	\end{equation}
	Since the condition (\ref{sec1:gam_bd}) holds for all $l\geq M_k$, one gets $\bar{\gamma}_{M_k|k}\geq\underline{\gamma}_{M_k|k}$. Hence, $\bar{\gamma}_{M_k|k+1}\geq\underline{\gamma}_{M_k|k+1}$ by following (\ref{sec1:gamma_lemma1_1}). In addition, from the condition $\mathcal{Z}_{l+1|k}\subseteq\mathcal{Z}_{l|k}$, we have $\bar{\gamma}_{l+1|k}\geq \bar{\gamma}_{l|k}$ and $\underline{\gamma}_{l+1|k}\leq \underline{\gamma}_{l|k}$. Therefore, there must exist a non-empty set 
	\begin{equation*}
		\mathbb{M}_{k+1} = \{m\in\mathbb{N}_0^{M_k}|\bar{\gamma}_{M_k-m|k+1}\geq\underline{\gamma}_{M_k-m|k+1}\}.
	\end{equation*}
	According to Proposition \ref{sec1:bar_alpha}, $M_{k+1}$ can be chosen as $M_{k+1} = M_k-m$. Therefore, we have $M_{k+1}\leq M_k$. \QEDA
\end{pf*}
\begin{rem}
	From Proposition \ref{sec1:bar_alpha}, it can be seen that extra $M_k$ {steps are} required to steer $\alpha_{l|k}$ into the terminal set $\mathcal{A}_k$. Hence, the prediction horizon is extended from $N$ to $N+M_k$. Based on Lemmas \ref{sec1:lemma:uncer_set} and \ref{sec1:horizon_mode2}, it can be derived that the sequence $\{ M_k\}$ is non-increasing. Hence, when $k$ increases, the computational complexity of MPC optimization problem is non-increasing. 
\end{rem}

To find the terminal set for the nominal state $z_{N|k}$, we have the following assumption:
\begin{assum}\label{sec1:assumption:RPI_set}
	Let $\mathcal{Z}_k$ and $\mathcal{Z}_{k+1}$ denote the MRPI sets with respect to the uncertainty set $\Theta_{k+1}$ and $\Theta_{k+2}$, respectively. Then the following condition holds
	\begin{equation}\label{sec1:P_assum}
		\phi(\theta) x\in\mathcal{Z}_{k+1}, \forall(x, \theta)\in\mathcal{Z}_{k}\times \Theta_{k+2}
	\end{equation}
	if $\Theta_{k+2}\subseteq \Theta_{k+1}$.
\end{assum}
\begin{rem}
	To compute the set $\mathcal{Z}_{k+1}$ satisfying the condition (\ref{sec1:P_assum}), we can compute the RPI set  $\bar{\mathcal{Z}}_{k+1}$ 
	by following Algorithm 1 in \cite{pluymers_efficient_2005ACC} without considering (\ref{sec1:P_assum}). Then starting with $\bar{\mathcal{Z}}_{k+1}$, $\mathcal{Z}_{k+1}$ can be computed by solving the linear programming problem with the additional constraint (\ref{sec1:P_assum}). In addition, given $\mathcal{Z}_k,\Theta_{k+1}$ and $\Theta_{k+2}$ with $\Theta_{k+2}\subseteq\Theta_{k+1}$, there always exists one $\mathcal{Z}_{k+1}$ such that (\ref{sec1:P_assum}) holds. A simple example is to choose $\mathcal{Z}_{k+1}$ as $\mathcal{Z}_{k+1} = \mathcal{Z}_k$ directly.
\end{rem}

\begin{assum}\label{sec1:assum:horizon2}
	Let $M_k,M_{k+1},\mathcal{A}_{k}$ and $\mathcal{A}_{k+1}$ are the horizons and invariant sets satisfying Proposition \ref{sec1:bar_alpha} with respect to uncertainty sets $\Theta_{k+1}$ and $\Theta_{k+2}$, respectively. Given $M_k$ and ${A}_{k}$, if the condition $\Theta_{k+2}\subseteq\Theta_{k+1}$ holds, there exist $M_{k+1}$ and $\mathcal{A}_{k+1}$ such that $M_k\geq M_{k+1}$ and $\mathcal{A}_{k}\subseteq\mathcal{A}_{k+1}$.
\end{assum}

{According to Proposition \ref{sec1:bar_alpha}, the feasible solution set of $\gamma_{k}$ in (\ref{sec1:gam_bd}) becomes larger when $l$ increases. Therefore, the larger invariant set $\mathcal{A}_k$ can be found by choosing the larger horizon $\mathcal{M}_k$. In addition, it follows from (\ref{sec1:uncer_set}) that $\Theta_{k+2}\subseteq\Theta_{k+1}$ for all $k\geq0$.  Let $\gamma_k = \bar{\gamma}_{M_k|k}, \gamma_{k+1} = \bar{\gamma}_{M_{k+1}|k+1}$ and $M_{k+1} = M_k$, then we have $\mathcal{A}_{k}\subseteq\mathcal{A}_{k+1}$ since ${\gamma}_{k+1}\geq {\gamma}_{k}$. Therefore, given $M_{k}$ and $\mathcal{A}_{k}$, we can always find $M_{k+1}$ and $\mathcal{A}_{k+1}$ such that Assumption \ref{sec1:assum:horizon2} holds. As a result, the computational complexity of MPC optimization problem is still non-increasing under this assumption.}

Suppose that the RPI set $\mathcal{Z}_{k}$ has the polyhedral form $\mathcal{Z}_k = \{x| V_kx\leq\mathbf{1}\}$, then the terminal constraints for the systems in (\ref{sec1:sys_mode2}) are summarized as follows:
\begin{subequations}\label{sec1:z_terminal}
	\begin{align}
		&V_{k} z_{N|k} + D_{k}\alpha_{N|k}\leq \mathbf{1}\label{sec1:z_terminal_var_1},\\
		&0\leq \alpha_{N+M_k|k}\leq \gamma_k \mathbf{1},\label{sec1:z_terminal_var_2}
	\end{align}
\end{subequations}
where $D_{k}$ is a non-negative matrix satisfying $D_{k}V = V_{k}$.
%%%%%%%%%%%%%%%%%%%%%%%%%%%%%%%%%%%%%%%%%%%%%%%%%%%%%%%%%%%%%%%%%%%%%%%%%%%%%%%%%%%%%%%%%%%%%%%%%%%%%%%%%%%%%%%%%%%%%%%%%%%%%%%%%%%%%%%%
\subsection{Construction of the cost function}
Let $\textbf{v}_k = \colt(v_{0|k},v_{1|k},v_{2|k},\cdots,v_{N-1|k})$. Define $E$ and $T$ as shift matrices such that $v_{0|k}=E\textbf{v}_k$ and $\textbf{v}_{k+1} = T\textbf{v}_{k}$, then the prediction of $z_{l|k}$ can be written as $\xi_{l+1|k}=\Psi_{k+1}\xi_{l|k}$, where $\xi_{l|k} = \begin{bmatrix}z_{l|k}\\\textbf{v}_k\end{bmatrix}, \Psi_{k+1} = \begin{bmatrix}\phi_{k+1}&B_{k+1}E\\0&T\end{bmatrix}.$ Similarly, the real system state $x_{l|k}$ can be predicted by using the following dynamics $\bar{\xi}_{l+1|k}=\Psi^*\bar{\xi}_{l|k}$,  where $\bar{\xi}_{l|k} = \begin{bmatrix}x_{l|k}\\\textbf{v}_k\end{bmatrix}$ and $\Psi^* = \begin{bmatrix}\phi^*&B^*E\\0&T\end{bmatrix}.$ In this work, the objective is to minimize a cost function with a quadratic form $\bar{J}_k= \sum_{i=0}^{\infty}(x_{i|k}^\T Qx_{i|k}+u_{i|k}^\T Ru_{i|k})$, where $Q>0$ and $R>0$ are penalty matrices for the state and input, respectively. Note that the cost function $\bar{J}_k$ can be equivalently represented by $\bar{J}_k=\xi_{0|k}^\T W^*\xi_{0|k}$ where $W^*$ is the solution of a Lyapunov equation
\begin{equation}\label{sec1:lyp_real}
	\begin{array}{l}
		(\Psi^*)^\T W^*(\Psi^*)-W^*+\bar{Q}=0
	\end{array}
\end{equation}  with $
\bar{Q} = \begin{bmatrix}
Q + K^\T R K & K^\T R E\\E^\T R K &E^\T R E \end{bmatrix}$.  Since $\phi^*$ is unknown, we cannot find the matrix $W^*$ exactly. Alternatively, we consider an over approximation of $\bar{J}_k$ based on the uncertainty set updated at each time instant. 
\begin{lem}\label{sec1:lemma:Wk}
	Define a new cost function $J_k$ as $J_k = \xi_{0|k}^\T W_{k+1}\xi_{0|k}$, where $W_{k+1}$ is a positive definite matrix, then $J_k \geq \bar{J}_k$ if the following condition
	\begin{equation}\label{sec1:wk_lyp_geq}
		W_{k+1} \geq \begin{bmatrix}\phi(\theta)&B(\theta)E\\0&T\end{bmatrix}^\T W_{k+1} \begin{bmatrix}\phi(\theta)&B(\theta)E\\0&T\end{bmatrix}+\bar{Q}
	\end{equation}
	holds for all $\theta\in\Theta_{k+1}$.
\end{lem}
\begin{pf*}{Proof}
	From Lemma \ref{sec1:lemma:uncer_set}, we have $\theta^*\in\Theta_{k+1}$. Then following (\ref{sec1:wk_lyp_geq}) yields $W_{k+1}\geq (\Psi^*)^\T W_{k+1}(\Psi^*)+\bar{Q}$.  By substituting $\bar{Q} = W^* - (\Psi^*)^\T W(\Psi^*)$ into the above equation, we have $W_{k+1}-W^*\geq (\Psi^*)^\T (W_{k+1}-W^*)(\Psi^*)\geq 0$. In addition, $J_k-\bar{J}_k = \xi_{0|k}^\T W_{k+1}\xi_{0|k} - \bar{\xi}_{0|k}^\T W^*\bar{\xi}_{0|k}$. Since $\bar{\xi}_{0|k} = {\xi}_{0|k}$ and $W_{k+1}-W^* \geq 0$, it can {be concluded} that $J_k\geq\bar{J}_k$ for all $\theta\in\Theta_{k+1}$. 
\QEDA\end{pf*}
%For the sequence of weighting matrices $W_k$, we have the following assumption.
\begin{assum}\label{sec1:assumption:cost_weight}
	Let $W_{k+1}$ denote the weighting matrix at time $k$, if $\Theta_{k+1}\subseteq\Theta_{k}$, then the following condition holds for all $k\geq 0$  
	\begin{equation}
		\xi_{0|k}^\T {W_{k+1}}\xi_{0|k}\leq\xi_{0|k}^\T {W_{k}}\xi_{0|k}.\label{sec1:stability_condition}
	\end{equation}
\end{assum}
\begin{rem}
	{Following (\ref{sec1:est_w})-(\ref{sec1:uncer_set}), it can be guaranteed that  $\Theta_{k+1}\subseteq\Theta_{k}$ for all $k\geq 0$. Given $W_k$,by imposing (\ref{sec1:stability_condition}) as an additional constraint for the {\normalfont\text{LMI}} problem used for computing $W_{k+1}$, we can find a $W_{k+1}$ satisfying the condition (\ref{sec1:stability_condition}).} An example of formulating the {\normalfont\text{LMI}} problem can be found in \cite{kouvaritakis2015model} for details.
\end{rem}
%%%%%%%%%%%%%%%%%%%%%%%%%%%%%%%%%%%%%%%%%%%%%%%%%%%%%%%%%%%%%%%%%%%%%%%%%%%%%%%%%%%%%%%%%%%%%%%%%%%%%%%%%%%%%%%%%%%%%%%%%%%%%%%%%%%%%%%%%
\subsection{Adaptive MPC algorithm}
According to the developed terminal sets and cost function, the adaptive MPC algorithm is based on the following MPC optimization problem:
\begin{subequations}
	\begin{align*}
		\mathbb{P}:	\min_{\textbf{v}_k,\{\alpha_{l|k}\}}~~&J_k = \xi_{0|k}^\T {W_{k+1}}\xi_{0|k}\\
		\text{s.t.}~~ &z_{0|k}=x_k\\
		&(\ref{sec1:u_dual}),(\ref{sec1:sys_nol}),(\ref{sec1:tube1}),(\ref{sec1:tube2}),(\ref{sec1:z_terminal_var_1}),(\ref{sec1:z_terminal_var_2})
	\end{align*}
\end{subequations}

At time instant $k$, we update the estimation of the unknown parameters and the uncertainty set based on new measurements, then reformulate the optimization problem $\mathbb{P}$. Note that the reformulation of $\mathbb{P}$ with respect to the new estimation is not necessary if the estimation error is sufficiently small. To reduce redundant estimating actions, we introduce a termination criterion for the proposed estimator. Let $\epsilon_x>0$ and $\epsilon_r>0$ denote the tolerances for the state estimation error and the error bound of parameter estimation, then the proposed adaptive MPC algorithm is summarized in Algorithm \ref{sec1:alg:adaptiveMPC}. 

\begin{algorithm}[!t]
	\caption{The Adaptive MPC algorithm}
 	\begin{algorithmic}[1]\label{sec1:alg:adaptiveMPC}
%		\REQUIRE At time $k$, given $x_k,\Theta_{k},\hat{\theta}_{k},\mathcal{Z}_{k-1},M_{k-1}, \mathcal{A}_{k-1}$ and $W_k$. 
		\REQUIRE Given initial conditions $x_0,\Theta_{0}$ and weighting matrices $Q,R$, determine the prestabilizing feedback gain $K$ and MRPI set $\mathcal{Z}_{0}$. Compute the terminal set $\mathcal{A}_{0}$ and the horizon $M_{0}$ according to Proposition \ref{sec1:bar_alpha}. Calculate the weighting matrix $W_0$ satisfying (\ref{sec1:wk_lyp_geq}).
		%\REQUIRE Given initial conditions $x_0,\Theta_{0}$, choose weighting matrices $Q,R$. Determine the prestabilizing feedback gain $K$. Calculate the weighting matrix $W_0$ satisfying (\ref{sec1:wk_lyp_geq}).
		\FOR {each time instant $k=0,1,2,\cdots$} 
		\IF {$\|\tilde{x}_k\|\geq \epsilon_x$ or $\mathcal{V}_{k}\geq \epsilon_r$}
		\STATE Calculate ${\hat{\theta}}_{k+1}$ and $\Theta_{k+1}$ by using (\ref{sec1:est_w})-(\ref{sec1:uncer_set}).\;
		\STATE Compute $M_{k}, \mathcal{A}_{k}, W_{k+1}$ and $\mathcal{Z}_{k}$ with respect to $\Theta_{k+1}$ such that Assumptions \ref{sec1:assumption:RPI_set}, \ref{sec1:assum:horizon2} and \ref{sec1:assumption:cost_weight} hold.
		\ELSE
		\STATE Let ${\hat{\theta}}_{k+1} = {\hat{\theta}}_{k},\Theta_{k+1} = \Theta_k,\mathcal{Z}_{k} = \mathcal{Z}_{k-1},W_{k+1} = W_k, M_{k} = M_{k-1}$ and $ \mathcal{A}_{k} = \mathcal{A}_{k-1}$.
		\ENDIF
		\STATE Reformulate and solve the optimization problem $\mathbb{P}$ based on $\hat{\theta}_{k+1}$ and $\Theta_{k+1}$ to obtain $\textbf{v}^*_k$.\;, 
		\STATE Calculate the control input as $u_k = Kx_k+v^*_{0|k}$, and then implement $u_k$ to the system.\;
		\ENDFOR
	\end{algorithmic}
\end{algorithm} 

\begin{thm}\label{sec1:theorem:feasbility}
	Suppose that Assumptions \ref{sec1:assumption:RPI_set}, \ref{sec1:assum:horizon2} and \ref{sec1:assumption:cost_weight} hold, and there is a feasible solution to the optimal control problem $\mathbb{P}$ when $k = 0$. Then $\mathbb{P}$ is recursively feasible by following Algorithm \ref{sec1:alg:adaptiveMPC}.
\end{thm}

\begin{pf*}{Proof}
	%Since $\hat{\theta}_k$ is updated by following (\ref{sec1:est_w})-(\ref{sec1:uncer_set}), from Lemma \ref{sec1:lemma:uncer_set}, we have $\Theta_{k+1}\subseteq\Theta_{k}$.
	Suppose that $\mathbb{P}$ is feasible at time $k$. Let $\textbf{v}_k^*$ and $\boldsymbol{\alpha}^*_k = \{\alpha^*_{l|k}\}_{l\in\mathbb{N}_0^{N+M_k}}$ denote the optimal solution of the MPC problem at time $k$. {$\{z^*_{l|k},\mathcal{S}^*_{l|k} =\{ e_{l|k}| V_ke_{l|k}\leq\alpha^*_{l|k}\},\mathcal{X}^*_{l|k} = z^*_{l|k}\oplus S^*_{l|k}\}_{l\in\mathbb{N}_0^{N+M_k}}$ are the corresponding nominal states, error tubes and state tubes, respectively.} Define a candidate input sequence at time $k+1$ as $\bar{\textbf{v}}_{k+1} = \{{v}_{1|k}^*,{v}_{2|k}^*,\cdots,{v}_{N-1|k}^*,0\}$. 
	
	Two cases are investigated to prove this theorem.\\
	\textbf{Case (1)}: Suppose that the estimation termination criterion in Algorithm \ref{sec1:alg:adaptiveMPC} is not satisfied. Based on $z_{0|l+1}$ and $\bar{v}_{l|k+1}$, we firstly construct the following sequence $\bar{\boldsymbol{\alpha}}_{k+1}=\{\alpha_{l|k+1}\}_{l\in\mathbb{N}_0^{N+M_{k+1}-1}}$ such that $\mathcal{X}_{l|k+1} = \mathcal{X}^*_{l+1|k}$. Let $\alpha_{N+M_{k+1}|k+1} = \underset{j\in\mathbb{N}_0^{n_c}}{\max}\{H_{k+2}^j\alpha_{N+M_{k+1}-1|k+1}+V\Delta\phi_{k+2}^jz_{N+M_{k+1}-1|k+1}\}$, we show that $\{\bar{\textbf{v}}_{k+1},\bar{\boldsymbol{\alpha}}_{k+1}\}$ is a feasible solution for $\mathbb{P}$ in the following.
	\begin{itemize}
		\item For $l\in\mathbb{N}_0^{N+M_{k+1}-1}$, since $\mathcal{X}_{l|k+1} = \mathcal{X}^*_{l+1|k}$, we have $\{z^*_{l+1|k},S^*_{l+1|k}\}$ satisfying the condition $z^*_{l+1|k}\oplus S^*_{l+1|k} = z_{l|k+1}\oplus S_{l|k+1}$, which verifies that the candidate sequence $\{z_{l|k+1},\alpha_{l|k+1},\bar{v}_{l|k+1}\}_{l\in\mathbb{N}_0^{N+M_{k+1}-1}}$ satisfies the constraints (\ref{sec1:tube1}) and  (\ref{sec1:tube2}).
		
		\item When $l=N$, it follows form (\ref{sec1:z_terminal_var_1}) that $V_{k} z^*_{N|k} + D_{k}Ve_{N|k}\leq \mathbf{1}$. By using $D_kV = V_k$, we have $V_{k}(z^*_{N|k}+e_{N|k}) = V_{k}x_{N|k} \leq 1$, implying that $\mathcal{X}^*_{N|k}\subseteq\mathcal{Z}_{k}$. As aforementioned, $\mathcal{X}_{N-1|k+1} =\mathcal{X}^*_{N|k}$, then $\mathcal{X}_{N-1|k+1}\subseteq\mathcal{Z}_{k}$. Since $\mathcal{Z}_{k}$ is an RPI set, $\bar{v}_{N-1|k+1} = 0$ and $\Theta_{k+2}\subseteq\Theta_{k+1}$, it yields that $\mathcal{X}_{N|k+1}\subseteq \phi_{k+2}\mathcal{Z}_{k}\subseteq\mathcal{Z}_{k+1}$ by following Assumption \ref{sec1:assumption:RPI_set}. Hence, we have $z_{N|k+1}+e_{N|k+1}\in\mathcal{Z}_{k+1}$ for all admissible $e_{N|k+1}$. As a result, the constraint (\ref{sec1:z_terminal_var_1}) is satisfied.
		
		\item When $ l = N+M_{k+1}$, taking the infinity norm of (\ref{sec1:alpha_dyn}) we have $\|\alpha_{l|k+1}\|_\infty\leq \underset{j\in\mathbb{N}_0^{n_c}}{\max} \{ \|H_{k+2}^j\|_\infty\|\alpha_{l-1|k+1}\|_\infty+\|\bar{c}^j_{l-1|k+1}\|_\infty\}$. Since $\mathcal{X}_{l-1|k+1} = \mathcal{X}_{l|k}^*$ and $\Theta_{k+2}\subseteq\Theta_{k+1}$, following Proposition \ref{sec1:bar_alpha} and Assumption \ref{sec1:assum:horizon2}, it is concluded that  $\|\alpha_{l|k+1}\|_\infty\leq \gamma_{k}\leq\gamma_{k+1}$. Hence, the constraint (\ref{sec1:z_terminal_var_2}) is satisfied.
	\end{itemize}
	\textbf{Case (2)}: Suppose that the estimation termination criterion in Algorithm \ref{sec1:alg:adaptiveMPC} is satisfied. Then we have 
	 ${\hat{\theta}}_{k+2} = {\hat{\theta}}_{k+1},~\mathcal{Z}_{k+1} = \mathcal{Z}_k, W_{k+2} = W_{k+1}, \gamma_{k+1} = \gamma_k$ and $M_{k+1} = M_k$. The recursive feasibility can be proved by constructing the following candidate sequence $\bar{\textbf{v}}_{k+1}, \{\alpha_{1|k}^*,\alpha_{2|k}^*,\cdots,\alpha_{N+M_k|k}^*,\underset{j\in\mathbb{N}_0^{n_c}}{\max}\{H_{k+1}^j\alpha_{N+M_k|k}+\\V\Delta\phi_{k+1}^jz_{N+M_k|k}\}\}$.

	In summary, there is a feasible solution for the optimal control problem $\mathbb{P}$ at time $k+1$ if it is feasible at time $k$. Therefore $\mathbb{P}$ is proved to be recursively feasible.
\QEDA\end{pf*}

\begin{thm}\label{sec1:theorem:stability}
	Suppose that Assumptions \ref{sec1:assumption:RPI_set}, \ref{sec1:assum:horizon2} and \ref{sec1:assumption:cost_weight} hold, then the system in (\ref{sec1:sys}) {in closed-loop} is asymptotically stable by applying the adaptive MPC Algorithm \ref{sec1:alg:adaptiveMPC}.
\end{thm}

\begin{pf*}{Proof}
	{To prove this theorem, in the following, we show that the optimal cost $J^*_k$ is a Lyapunov function for the system in (\ref{sec1:sys}) in closed-loop with Algorithm \ref{sec1:alg:adaptiveMPC}.}
		
	{\textbf{Case (1)}: Suppose that the estimation termination criterion in Algorithm \ref{sec1:alg:adaptiveMPC} is not satisfied. Let $z_{0|k+1} = x_{k+1}, \xi_{0|k+1} =\colt(z_{0|k+1},\bar{\textbf{v}}_{k+1}), \xi_{0|k} = \colt(z_{0|k}^*,{\textbf{v}}_{k}^*)$ and $\bar{J}_{k+1} = \xi_{0|k+1}^\T {W_{k+2}}\xi_{0|k+1}$, based on Lemma \ref{sec1:lemma:Wk}, we have
	\begin{equation*}
		\begin{array}{l}
			~~~{\xi}_{0|k+1}^\T W_{k+1}{\xi}_{0|k+1} - J_k^*\\
			= \xi_{0|k}^\T(\Psi^*)^\T W_{k+1} \Psi^*\xi_{0|k} - \xi_{0|k}^\T W_{k+1}\xi_{0|k}\\
			\leq -\xi_{0|k}^\T\bar{Q}\xi_{0|k}\\
			 = -z_{0|k}^\T{Q}z_{0|k} - u_{0|k}^\T{R}u_{0|k}			
		\end{array}
	\end{equation*}
	Since $Q$ and $R$ are positive definite and $z_{0|k} = x_k$, it can be derived that $\bar{\xi}_{0|k+1}^\T W_{k+1}\bar{\xi}_{0|k+1} - J_k^*\leq-x_{k}^\T{Q}x_{k}- u_{0|k}^\T{R}u_{0|k}$. In addition, from Assumption \ref{sec1:assumption:RPI_set}, we have $\bar{J}_{k+1}=\bar{\xi}_{0|k+1}^\T {W_{k+2}}\bar{\xi}_{0|k+1} \leq \bar{\xi}_{0|k+1}^\T W_{k+1}\bar{\xi}_{0|k+1}$, which yields $J_{k+1}^* - J_k^*\leq \bar{J}_{k+1} - J_k^* \leq -x_{k}^\T{Q}x_{k}-u_{k}^\T{R}u_{k}\leq 0,~\forall x_k\neq0$ and $u_k\neq 0$. Since $W_k$ is positive definite, $J_k^*$ is a Lyapunov function for the system in (\ref{sec1:sys}).} 
	
{	\textbf{Case (2)}: Suppose that the estimation termination criterion in Algorithm 1 is satisfied. Then we have
		${\hat{\theta}}_{k+2} = {\hat{\theta}}_{k+1},~\mathcal{Z}_{k+1} = \mathcal{Z}_k, W_{k+2} = W_{k+1}, \gamma_{k+1} = \gamma_k$ and $M_{k+1} = M_k$. By repeating the above procedure, we can prove that $J_k^*$ is a Lyapunov function.}
	
{In summary, the optimal cost function $J_k^*$ is a Lyapunov function for the system in (\ref{sec1:sys}) in closed-loop with Algorithm \ref{sec1:alg:adaptiveMPC}. Hence, the closed-loop system is asymptotically stable.} 
\QEDA\end{pf*}
\begin{rem}
{	Note that, unlike the robust method in \cite{fleming_robust_2015TAC}, the propagation of homothetic tube $\mathcal{S}_{l|k}$ (\ref{sec1:tube}) in our proposed method depends on the estimation $\hat{\theta}_{k+1}$ and $\hat{\Theta}_{k+1}$. In addition, the nominal system in (\ref{sec1:sys_nol}), the terminal conditions in (\ref{sec1:z_terminal}) and  the weighting matrix $W_{k+1}$ are also updated based on the estimation of uncertainty at each time instant. By following (\ref{sec1:est_w})-(\ref{sec1:uncer_set}), the non-increasing properties on the proposed estimation scheme are guaranteed. Therefore, the resulting adaptive MPC scheme can reduce conservatism compared with the original robust MPC method. The numerical simulations will elaborate this argument.}     
\end{rem}

\begin{rem}
	As shown in Algorithm \ref{sec1:alg:adaptiveMPC}, when updating the parameter estimate $\hat{\theta}_k$ and uncertainty set $\Theta_{k}$, we need to re-compute $M_k,\mathcal{Z}_{k+1}$ and $W_{k+1}$, which is relatively computationally expensive. For some problems which have the strict requirement on the computational load, a solution to reduce the computational complexity is to choose the relatively large $\epsilon_x$ and $\epsilon_r$. An alternative is to omit the update of terminal conditions and cost function by setting $M_k = M_0,\mathcal{Z}_{k} = \mathcal{Z}_{0}$ and $W_{k} = W_0$ for all $k\geq 0$. Due to the fact that $\Theta_{k+1}\subseteq\Theta_{k}\subseteq\Theta_{0}$,	this strategy can significantly reduce the computational {load} with guaranteed closed-loop stability and recursive feasibility, but results in a relatively conservative control performance. {Note that the recursive updates of system model and uncertainty set are considered in the tube propagation, and thus, this simplified method still has less conservative closed-loop performance compared with the robust MPC method.} The numerical simulation will demonstrate this argument.
\end{rem}
%%%%%%%%%%%%%%%%%%%%%%%%%%%%%%%%%%%%%%%%%%%%%%%%%%%%%%%%%%%%%%%%%%%%%%%%%%%%%%%%
%%%%%%%%%%%%%%%%%%%%%%%%%%%%%%%%%%%%%%%%%%%%%%%%%%%%%%%%%%%%%%%%%%%%%%%%%%%%%%%%
\section{Simulation Results}\label{sec1:simulation}
In this section, a numerical example is presented to show the effectiveness of proposed adaptive MPC algorithms. The numerical test is conducted in Matlab, where the MPC optimization problem is formulated and solved by using Yalmip \cite{lofberg_yalmip_2004CACSD}.

We consider the following example for testing:
\begin{equation*}
\begin{array}{l}
A_0 = \begin{bmatrix}0.42&-0.28\\0.02&0.6\end{bmatrix},
A_1 = \begin{bmatrix}-0.12&-0.08\\-0.12&-0.17\end{bmatrix},
A_2 = -A_1,\\
B_0 = \begin{bmatrix}0.3 & -0.4\end{bmatrix}^\T,
B_1 = \begin{bmatrix}0.04&-0.08\end{bmatrix}^\T, 
B_2 = -1.5B_1.
\end{array}
\end{equation*}
%The initial uncertainty set is given by $\Theta_{0}= \{\theta\in\mathbb{R}^2 | \|\theta\|\leq 1\}$. The system is subject to input and state constraints $\{x|\|x\|_{\infty}\leq 17 \}$ and  $\{ u| \|u\|_\infty\leq 4\}$. 
$\Theta_{0}= \{\theta\in\mathbb{R}^2 | \|\theta\|\leq 1\}, \{x|\|x\|_{\infty}\leq 17 \}$ and $\{ u| \|u\|_\infty\leq 4\}$. The weighting matrices are chosen as $Q = I_2$ and $R=1$. By following \cite{kouvaritakis2015model}, the prestabilizing feedback gain is chosen as $K = [-0.4187~1.1562]$. Set the prediction horizon $N = 10$, then the horizon and terminal region are derived as $ M_0 = 3$ and $\gamma_0 = 0.4266$. The parameters used in Algorithm \ref{sec1:alg:adaptiveMPC} are given in the following $\epsilon_r=0.001,\epsilon_x = 0.001,\lambda = 0.5$ and $\Gamma_0=0.15I_2$. 

\begin{figure}
	\centering
	\includegraphics[width=0.8\linewidth]{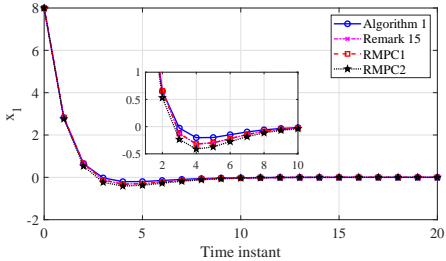}
	\caption{The time evolution of the system state $x_1$.}
	\label{sec1:fig:trajx1comv2}
\end{figure}

\begin{figure}
	\centering
	\includegraphics[width=0.8\linewidth]{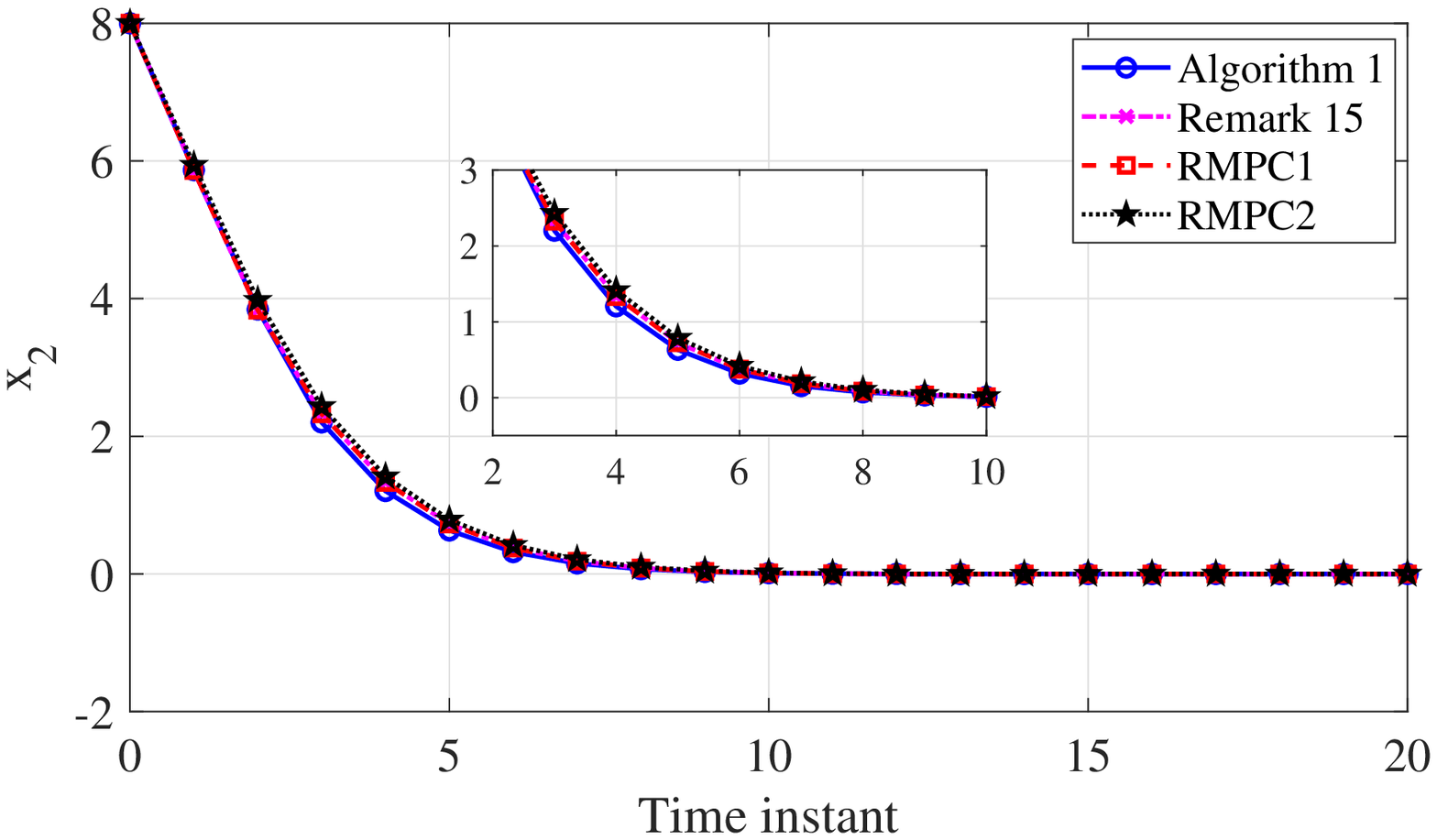}
	\caption{The time evolution of the system state $x_2$.}
	\label{sec1:fig:trajx2comv2}
\end{figure}

\begin{figure}
	\centering
	\includegraphics[width=0.7\linewidth]{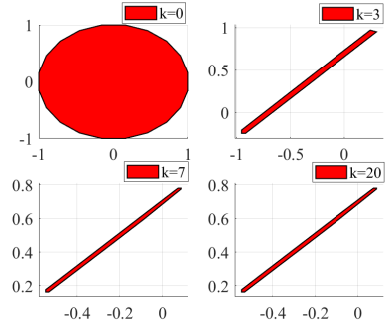}
	\caption{The estimated uncertainty set $\Theta$ obtained at $k = 0, 3,7,20$.}
	\label{sec1:fig:setcomv2}
\end{figure}

The robust MPC method in \cite{fleming_robust_2015TAC} (RMPC1) and \cite{lorenzen_robust_2019Automatica} (RMPC2) are introduced for the purpose of comparison. The initial point is set as $x_0 = [8,8]^\T$. The real system parameter $\theta^* = [-0.2,0.5]^\T$ is given to evaluate the proposed parameter estimator. Figs. \ref{sec1:fig:trajx1comv2} and \ref{sec1:fig:trajx2comv2} show the trajectories of system state and control input obtained by {applying different control} methods. From these figures, it can be seen that the recursive feasibility can be guaranteed by using these methods while the proposed method can accelerate the convergence of system state. To further compare the control performances of different MPC formulations, we introduce the following  index $\bar{J}_p = \sum_{k=0}^{T_{stp}}(x_k^\T Qx_k+u_k^\T R u_k)/T_{stp}$, where $T_{stp}$ denotes the simulation time. The corresponding results are illustrated in Table \ref{sec1:tab:averagecost}, implying that the proposed method can achieve the less conservative performance. The polytopic approximation of uncertainty sets obtained at time $k = 0, 3,7,20$ are depicted in Fig. \ref{sec1:fig:setcomv2}.  It can be seen that the estimate of uncertainty set is non-increasing, and finally converges to a fixed set, which verifies the proposed results.

\begin{figure}
	\centering
	\includegraphics[width=0.75\linewidth]{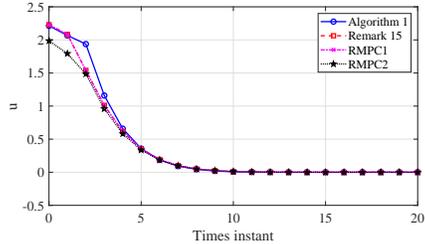}
	\caption{Trajectories of control input $u$.}
	\label{sec1:fig:trajucomv2}
\end{figure}

\begin{table}[htbp]
	\centering
	\begin{tabular}{@{}ccccc@{}}
		\toprule
		& Algorithm \ref{sec1:alg:adaptiveMPC} & Remark 15 &RMPC1 & RMPC2 \\ \midrule
		$\bar{J}_p$ & 9.2023 & 9.2524 & 9.2524 & 9.3747 \\ \bottomrule                               
	\end{tabular}
	\caption{The comparison of system performance.}
	\label{sec1:tab:averagecost}
\end{table}
%%%%%%%%%%%%%%%%%%%%%%%%%%%%%%%%%%%%%%%%%%%%%%%%%%%%%%%%%%%%%%%%%%%%%%%%%%%%%%%%%%%%%%%%%%%%%%%%%%%%%%%%%%%%%%%%%%%%%%%%%%%%
\section{Conclusion}\label{sec1:conclusion}
In this paper, we have investigated adaptive MPC for constrained linear systems subject to multiplicative uncertainties. An online parameter estimator has been designed based on the RLS technique for simultaneous parameter identification and uncertainty set estimation. By integrating the proposed estimator with homothetic prediction tubes, the resulting tube-based adaptive MPC scheme is recursively feasible with recursive model updates, while giving rise to enhanced performance compared with the robust tube MPC method. The simplified version of the proposed adaptive MPC method was also given to provide a trade-off between conservatism and computational complexity. We have proven that the closed-loop system is asymptotically stable. Numerical simulations and comparison studies have been given to demonstrate the efficacy and advantages of the proposed adaptive MPC method. {On the other hand, the main limitation of proposed adaptive MPC approach comes from the polytopic over approximation of the uncertainty set employed in the construction of homothetic tubes, leading to an undesired increase in conservatism and computational complexity. Furthermore, this work considered the constant parametric uncertainties only, which potentially poses certain limitations to practical applications. The future research will focus on how to efficiently use poytopes for bounding the FSS of unknown parameters with a tight overestimation. Incorporating such an idea to develop adaptive MPC algorithms for handling time-varying multiplicative and additive disturbances is also an interesting direction for future research.}     

%%%%%%%%%%%%%%%%%%%%%%%%%%%%%%%%%%%%%%%%%%%%%%%%%%%%%%%%%%%%%%%%%%%%%%%%%%%%%%%%%
%\begin{ack}                               % Place acknowledgements
%Partially supported by the Roman Senate.  % here.
%\end{ack}

\bibliographystyle{plain}        % Include this if you use bibtex 
\bibliography{ref_sim}           % and a bib file to produce the 
                                 % bibliography (preferred). The
                                 % correct style is generated by
                                 % Elsevier at the time of printing.

%\appendix
%\section{A summary of Latin grammar}    % Each appendix must have a short title.
%\section{Some Latin vocabulary}         % Sections and subsections are supported  
                                        % in the appendices.
\end{document}